\documentclass{article}
\usepackage{graphicx} 
\usepackage[utf8]{inputenc}
\usepackage{titlesec}
\titleformat{\section}
{\normalfont\scshape}
{\thesection.}{0.5em}{\centering}
\titleformat{\subsection}
{\small\scshape}
{\thesubsection.}{0.5em}{\centering}
\titleformat{\subsubsection}
{\small\scshape}
{\thesubsubsection.}{0.5em}{\centering}
\usepackage{tikz-cd}
\usepackage{hyperref}
\hypersetup{
    colorlinks=true,
    linkcolor=black,    
    urlcolor=cyan,
    citecolor=black
    }
\usepackage{tikz}
\usepackage{comment}
\usepackage{amsfonts}
\usepackage{indentfirst}
\renewcommand{\abstract}{\small{\section*{\abstractname}}}
\usepackage{amsmath}
\usepackage{amssymb}
\usepackage{amsthm}
\usepackage{mathrsfs}
\newtheorem*{theorem*}{Theorem}
\newtheorem*{lemma*}{Lemma}

\newtheorem{definition}{Definition}

\usepackage{stmaryrd}
\usepackage{fancyhdr}
\fancypagestyle{plain}{%
  \fancyhf{} 
  \fancyfoot[C]{\tiny This a preprint of the work accepted for publication in \href{https://www.pleiades.online/}{Mathematical Notes}, \textcopyright\ 2024, the copyright holder indicated in the Journal.}
  
}
\newcommand{\hooklongrightarrow}{\lhook\joinrel\longrightarrow}

\usepackage[
backend=biber,
style=apa,
sorting=ynt
]{biblatex}
\title{\normalsize{\textbf{\uppercase{On $\mathfrak{L}_{\infty}$-liftings of the Gelfand-Naimark morphism}}}}
\author{\normalsize{\textsc{Vasily Melnikov}}}
\date{\normalsize{\textsc{March 2024}}}
\begin{document}
\maketitle
\begin{abstract}
     Let $M$ be a $W^{\ast}$-algebra acting on a separable complex Hilbert space $H$. We show that the inclusion of $M$ into $\mathscr{B}(H)$ factors through an $\mathfrak{L}_{\infty}$-space only if $M$ is a finite type $\mathrm{I}$ algebra.
\end{abstract}
\section{Introduction}\label{sec:intro}
For $C^{\ast}$-algebras, it is clear that if the Gelfand-Naimark morphism lifts by a $C^{\ast}$-morphism to an abelian $C^{\ast}$-algebra, then the underlying algebra is commutative. Thus, lifting properties of certain morphisms are equivalent to results about commutativity.
\par
This leads to the natural question of whether the existence of a lift \textit{in some sense} like above implies that the underlying algebra is `approximately' commutative. For example, suppose that the Gelfand-Naimark morphism lifts \textit{in a Banach-space theoretic sense} through an abelian $C^{\ast}$-algebra (for example, through a $C(K)$ space). Although this lift need not preserve the $C^{\ast}$-algebraic structure of the algebra, it is not unreasonable to suspect the Banach-space theoretic properties of $C(K)$ spaces (such as the Dunford-Pettis property) should carry over to $C^{\ast}$-algebras satisfying a lifting property of this form. This is not clear for arbitrary Banach spaces embedded into a $C(K)$ space; indeed, $C(K)$ does not have the hereditary Dunford-Pettis property if $K$ is not scattered (for example, if $K=[0,1]$).
\par
In this article, we show that a lifting property involving the $\mathfrak{L}_{\infty}$-spaces of Lindenstrauss-Pełczyński (which in particular include $C(K)$ spaces) for $W^{\ast}$-subalgebras of $\mathscr{B}(H)$, where $H$ is a separable complex Hilbert space, has several important implications. In particular, $\mathfrak{L}_{\infty}$-lifting properties are shown to be equivalent to certain other Banach-space theoretic condition on the von Neumann algebra, and carry important information on the structure of the underlying $W^{\ast}$-algebra. Indeed, from the main theorem of this note one obtains that the non-commutativity experienced by such $W^{\ast}$-algebras is finite-dimensional in nature.
\section{Main results}\label{sec:main}
Fix $p\in[1,\infty]$, and $\lambda\in[1,\infty)$. Recall that a Banach space $E$ is said to be an $\mathfrak{L}_{p,\lambda}$-space if the `local structure' of $E$ possesses essentially the same properties as that of $\ell^{p}$, up to the constant $\lambda$. The precise definition is the following.
\begin{definition}
    Let $p\in[1,\infty]$, and $\lambda\in[1,\infty)$. An infinite-dimensional Banach space $E$ is said to be an $\mathfrak{L}_{p,\lambda}$-space if, for any finite-dimensional subspace $F\subset E$, there is a finite-dimensional subspace $G\subset E$ containing $F$ such that
    \begin{equation*}
        d(G,\ell^{p}_{n})\leq\lambda,
    \end{equation*}
    for $n=\mathrm{dim}(G)$, where $d$ denotes the multiplicative Banach-Mazur distance, i.e.
    \begin{equation*}
        d(G,H)=\inf\left\{\Vert{T}\Vert\Vert{T^{-1}}\Vert:T:G\longrightarrow H\textrm{ is an isomorphism}\right\}.
    \end{equation*}
\end{definition}
We refer the reader to (\cite{bourg}) for additional details concerning the $\mathfrak{L}_{p,\lambda}$-spaces.
\par
Fix a complex Hilbert space $H$. Let $(M,\ell)$ be an operator space, where $\ell:M\longrightarrow\mathscr{B}(H)$ is the inclusion map.
\begin{definition}\label{defn:ac}
    Fix $\lambda\geq1$. The operator space $(M,\ell)$ is said to have the $\lambda$-Lindenstrauss-Pełczyński property if there exists an $\mathfrak{L}_{\infty,\lambda}$-space $E$ such that
    \[\begin{tikzcd}
	M && {\mathscr{B}(H)} \\
	& E
	\arrow["S", from=1-1, to=2-2]
	\arrow["T", from=2-2, to=1-3]
	\arrow["\ell", from=1-1, to=1-3]
\end{tikzcd}\]
for some pair $(S,T)$ of bounded linear operators $S:M\longrightarrow E$, $T:E\longrightarrow\mathscr{B}(H)$.
\end{definition}
Let us briefly explain the idea behind definition \ref{defn:ac}. Suppose that $M$ is a $C^{\ast}$-algebra acting on a complex Hilbert space $H$; denote by $\ell$ the natural inclusion map. The condition that $(M,\ell)$ has the $\lambda$-Lindenstrauss-Pełczyński property can be seen as an approximate substitute for commutativity. Indeed, if $E$ from definition \ref{defn:ac} is an $\mathfrak{L}_{\infty,\lambda}$-space, then $E''$ is isomorphic to $\ell^{\infty}(\kappa)$ for some cardinal $\kappa$ by (Proposition 1.33, \cite{bourg}) and Haydon's theorem (see \cite{hayd}). Thus, $E''$ is isomorphic to an abelian von Neumann algebra, so that the natural $C^{\ast}$-morphism\footnote{Given a $C^{\ast}$-algebra algebra $N$, $N''$ denotes the von Neumann enveloping algebra; equivalently, by the Sherman-Takeda theorem, $N''$ is the bidual of $N$ equipped with the Arens product.} $M\hooklongrightarrow\mathscr{B}(H)''$ factors (in a Banach-space theoretic sense) through an abelian von Neumann algebra. Of course, this lift need not preserve the $C^{\ast}$-algebraic structure of $M$.
\par
Actually, significantly more is true. Suppose that $M$ is a $W^{\ast}$-subalgebra of $\mathscr{B}(H)$ for some complex separable Hilbert space $H$; denote by $\ell$ the natural inclusion of $M$ into $\mathscr{B}(H)$. If $(M,\ell)$ has the $\lambda$-Lindenstrauss-Pełczyński property for some $\lambda\geq1$, then $(M,\ell)$ has the $\lambda$-Lindenstrauss-Pełczyński property for all $\lambda>1$, and one can choose for $E$ to be isometrically isomorphic to an abelian $W^{\ast}$-algebra. Furthermore, the $\lambda$-Lindenstrauss-Pełczyński property has deep consequences for the $W^{\ast}$-algebraic structure of $M$, as it makes the non-commutativity of $M$ finite-dimensional in nature. Indeed, one has the following theorem, which is the main result of this note.
\begin{theorem*}\label{thm:main}
    Let $M$ be a $W^{\ast}$-subalgebra of $\mathscr{B}(H)$ for some complex separable Hilbert space $H$, with $\ell$ denoting the inclusion map. Suppose that $(M,\ell)$ has the $\lambda$-Lindenstrauss-Pełczyński property for some $\lambda\geq1$. Then there exists a bounded sequence $\{n_{k}\}_{k}$ of nonnegative integers with
     \begin{equation}\label{eq:iso}
         M\simeq\bigoplus_{k=1}^{\infty}N_{k},
     \end{equation}
     where $N_{k}$ is a von Neumann algebra of type $\mathrm{I}_{n_{k}}$, and $M$ has the $\lambda$-Lindenstrauss-Pełczyński property for all $\lambda>1$. Furthermore, one may assume that $E$ from definition \ref{defn:ac} is an abelian $W^{\ast}$-algebra.
\end{theorem*}
The following lemma is simple, but important for the proof of the main theorem.
\begin{lemma*}\label{lem:heredit}
    Suppose that an operator space $(M,\ell)$ has the $\lambda$-Lindenstrauss-Pełczyński property. Then, for any subspace $N$ of $M$, $(N,\ell|_{N})$ has the $\lambda$-Lindenstrauss-Pełczyński property.
\end{lemma*}
\begin{proof}
    The following diagram commutes:
    \[\begin{tikzcd}
	N && {\mathscr{B}(H)} \\
	& E
	\arrow["{\ell|_{N}}", from=1-1, to=1-3]
	\arrow["{S|_{N}}"', from=1-1, to=2-2]
	\arrow["T"', from=2-2, to=1-3]
\end{tikzcd}\]
where $S$, $T$, and $E$ are from definition \ref{defn:ac}.
\end{proof}
The proof proceeds as follows. First, one may assume that $M$ is amenable. Then, using what are essentially Banach-space theoretical arguments, one may show that $M$ must have the Dunford-Pettis property, so that the results of (\cite{dunf}) apply to $M$. One then uses the above lemma and the universality of the hyperfinite type $\mathrm{II}_{1}$ factor to reduce the general case to the amenable case.
\begin{proof}[Proof of the main theorem.]
    Suppose that $(M,\ell)$ has the $\lambda$-Lindenstrauss-Pełczyński property for some $\lambda\geq1$. Let $E$, $S$, and $T$ be as in definition \ref{defn:ac}; note that $S$ induces an isomorphism of $M$ with a closed subspace of $E$.
    \par
    First, assume that $M$ is amenable. We claim that $M$ has the Dunford-Pettis property. Since $\mathfrak{L}_{\infty,\lambda}$-spaces have the Dunford-Pettis property (see Corollary 1.30, \cite{bourg}), it suffices to show that $M$ is isomorphic to a complemented subspace of $E$ (see Corollary 1.13, \cite{bourg}). By (Theorem 1, \cite{connamen}) and amenability, $M$ is complemented in $\mathscr{B}(H)$ by a projection $p$.  It suffices to show that $S$ induces an isomorphism of $M$ with a complemented subspace of $E$. The following diagram commutes.
    \[\begin{tikzcd}
	M & E & {\mathscr{B}(H)} \\
	M
	\arrow[Rightarrow, no head, from=1-1, to=2-1]
	\arrow["S", from=1-1, to=1-2]
	\arrow["T", from=1-2, to=1-3]
	\arrow["p", from=1-3, to=2-1]
\end{tikzcd}\]
    Thus, $S$ induces an isomorphism of $M$ with a complemented subspace of $E$, as desired. Since $M$ has the Dunford-Pettis property, it follows from (Theorem 3, \cite{dunf}) that (\ref{eq:iso}) holds.
    \par
    Assume now that $M$ is not necessarily amenable. We will show that $M$ is amenable. There exists a direct sum decomposition,
   \begin{equation*}
    M\simeq N\oplus L,
    \end{equation*}
    of $M$ into the direct sum of a type $\mathrm{I}$ algebra $N$, and a von Neumann algebra $L$ which decomposes into the direct sum of type $\mathrm{II}$ and type $\mathrm{III}$ algebras. Since $N$ is necessarily amenable, it suffices to show that $L\simeq\{0\}$. Suppose that $L$ is not trivial. Every type $\mathrm{II}$ or type $\mathrm{III}$ factor contains the hyperfinite type $\mathrm{II}_{1}$ factor $R$ as a subfactor. Thus, after applying the direct integral decomposition theorem of von Neumann and Murray, $L$ contains $A\otimes R$ for some nontrivial abelian von Neumann algebra $A$, so that $L$ contains $R$ as a $W^{\ast}$-subalgebra. Since $R$ is amenable and does not have the Dunford-Pettis property, it follows that $R$ cannot possess the $\lambda$-Lindenstrauss-Pełczyński property for any $\lambda\geq1$. This contradicts the lemma above. Thus,
    \begin{equation*}
    L\simeq\{0\},
    \end{equation*}
    as desired.
    \par
    Suppose that (\ref{eq:iso}) holds. There exists a sequence $\{L_{n}\}_{n}$ of hyperstonean topological spaces such that
    \begin{equation*}
        N_{k}\simeq C(L_{k})\otimes M_{n_{k}},
    \end{equation*}
    where $M_{n}=\mathscr{B}(\mathbb{C}^{n})$. Define 
    \begin{equation*}
        E=\bigoplus_{k=1}^{\infty}C(L_{k}\times\{1,\dots,n_{k}\}\times\{1,\dots,n_{k}\}).
    \end{equation*}
    Then $E$ can be equipped with the structure of an abelian von Neumann algebra, so that $E$ is isometrically isomorphic to $C(K)$ for some hyperstonean topological space $K$ by Gelfand duality. Clearly, $M$ is isomorphic as a Banach space to $E$, so that $M\hooklongrightarrow\mathscr{B}(H)$ factors through $E$. Fix $\lambda>1$. It is easy to see that $E$ is an $\mathfrak{L}_{\infty,\lambda}$-space, as it is a $C(K)$ space.
\end{proof}
\section*{Acknowledgments}
 The author would like to thank the anonymous referee for carefully reading the article, and providing several useful suggestions.
\printbibliography
\end{document}